
%

 \magnification=1200 \vsize=18cm \voffset=1cm \hoffset=1cm



\hsize=11.25cm
\parskip 0pt
\parindent=12pt

\catcode'32=9

\font\tenpc=cmcsc10
\font\eightpc=cmcsc8
\font\eightrm=cmr8
\font\eighti=cmmi8
\font\eightsy=cmsy8
\font\eightbf=cmbx8
\font\eighttt=cmtt8
\font\eightit=cmti8
\font\eightsl=cmsl8
\font\sixrm=cmr6
\font\sixi=cmmi6
\font\sixsy=cmsy6
\font\sixbf=cmbx6

\skewchar\eighti='177 \skewchar\sixi='177
\skewchar\eightsy='60 \skewchar\sixsy='60

\catcode`@=11

\def\tenpoint{%
  \textfont0=\tenrm \scriptfont0=\sevenrm \scriptscriptfont0=\fiverm
  \def\rm{\fam\z@\tenrm}%
  \textfont1=\teni \scriptfont1=\seveni \scriptscriptfont1=\fivei
  \def\oldstyle{\fam\@ne\teni}%
  \textfont2=\tensy \scriptfont2=\sevensy \scriptscriptfont2=\fivesy
  \textfont\itfam=\tenit
  \def\it{\fam\itfam\tenit}%
  \textfont\slfam=\tensl
  \def\sl{\fam\slfam\tensl}%
  \textfont\bffam=\tenbf \scriptfont\bffam=\sevenbf
  \scriptscriptfont\bffam=\fivebf
  \def\bf{\fam\bffam\tenbf}%
  \textfont\ttfam=\tentt
  \def\tt{\fam\ttfam\tentt}%
  \abovedisplayskip=12pt plus 3pt minus 9pt
  \abovedisplayshortskip=0pt plus 3pt
  \belowdisplayskip=12pt plus 3pt minus 9pt
  \belowdisplayshortskip=7pt plus 3pt minus 4pt
  \smallskipamount=3pt plus 1pt minus 1pt
  \medskipamount=6pt plus 2pt minus 2pt
  \bigskipamount=12pt plus 4pt minus 4pt
  \normalbaselineskip=12pt
  \setbox\strutbox=\hbox{\vrule height8.5pt depth3.5pt width0pt}%
  \let\bigf@ntpc=\tenrm \let\smallf@ntpc=\sevenrm
  \let\petcap=\tenpc
  \normalbaselines\rm}

\def\eightpoint{%
  \textfont0=\eightrm \scriptfont0=\sixrm \scriptscriptfont0=\fiverm
  \def\rm{\fam\z@\eightrm}%
  \textfont1=\eighti \scriptfont1=\sixi \scriptscriptfont1=\fivei
  \def\oldstyle{\fam\@ne\eighti}%
  \textfont2=\eightsy \scriptfont2=\sixsy \scriptscriptfont2=\fivesy
  \textfont\itfam=\eightit
  \def\it{\fam\itfam\eightit}%
  \textfont\slfam=\eightsl
  \def\sl{\fam\slfam\eightsl}%
  \textfont\bffam=\eightbf \scriptfont\bffam=\sixbf
  \scriptscriptfont\bffam=\fivebf
  \def\bf{\fam\bffam\eightbf}%
  \textfont\ttfam=\eighttt
  \def\tt{\fam\ttfam\eighttt}%
  \abovedisplayskip=9pt plus 2pt minus 6pt
  \abovedisplayshortskip=0pt plus 2pt
  \belowdisplayskip=9pt plus 2pt minus 6pt
  \belowdisplayshortskip=5pt plus 2pt minus 3pt
  \smallskipamount=2pt plus 1pt minus 1pt
  \medskipamount=4pt plus 2pt minus 1pt
  \bigskipamount=9pt plus 3pt minus 3pt
  \normalbaselineskip=9pt
  \setbox\strutbox=\hbox{\vrule height7pt depth2pt width0pt}%
  \let\bigf@ntpc=\eightrm \let\smallf@ntpc=\sixrm
  \let\petcap=\eightpc
  \normalbaselines\rm}
\catcode`@=12

\tenpoint
\font\tenbboard=msbm10
\def\bboard#1{\hbox{\tenbboard #1}}
\font\tengoth=eufm10

\catcode`\@=11
\def\pc#1#2|{{\bigf@ntpc #1\penalty \@MM\hskip\z@skip\smallf@ntpc%
    \uppercase{#2}}}
\catcode`\@=12

\def\pointir{\discretionary{.}{}{.\kern.35em---\kern.7em}\nobreak
   \hskip 0em plus .3em minus .4em }

\def\qed{\quad\raise -2pt\hbox{\vrule\vbox to 10pt{\hrule width 4pt
   \vfill\hrule}\vrule}}

\def\rem#1|{\par\medskip{{\it #1}\pointir}}

\def\vspace[#1]{\noalign{\vskip#1}}

\def\resume#1{\vbox{\eightpoint\narrower\narrower
\pc R\'ESUM\'E|.\quad #1}}
\def\abstract#1{\vbox{\eightpoint\narrower\narrower
\pc ABSTRACT|.\quad #1}}

\def\section#1{\goodbreak\par\vskip .3cm\centerline{\bf #1}
   \par\nobreak\vskip 3pt }

\long\def\th#1|#2\endth{\par\medbreak
   {\petcap #1\pointir}{\it #2}\par\medbreak}

\def\article#1|#2|#3|#4|#5|#6|#7|
    {{\leftskip=7mm\noindent
     \hangindent=2mm\hangafter=1
     \llap{{\tt [#1]}\hskip.35em}{\petcap#2}\pointir
     #3, {\sl #4}, {\bf #5} ({\oldstyle #6}),
     pp.\nobreak\ #7.\par}}

\def\livre#1|#2|#3|#4|
    {{\leftskip=7mm\noindent
    \hangindent=2mm\hangafter=1
    \llap{{\tt [#1]}\hskip.35em}{\petcap#2}\pointir
    {\sl #3}, #4.\par}}

\def\divers#1|#2|#3|
    {{\leftskip=7mm\noindent
    \hangindent=2mm\hangafter=1
     \llap{{\tt [#1]}\hskip.35em}{\petcap#2}\pointir
     #3.\par}}

\font\sevenbboard=msbm7
\def\setA{{\bboard A}}
\def\setB{{\bboard B}}
\def\setAA{\hbox{\sevenbboard A}}
\def\setZ{{\bboard Z}}
\def\setZq{{\bboard Z}[q,q^{-1}]}
\def\setZtq{{\bboard Z}[t,t^{-1},q,q^{-1}]}

\font\sevenbboard=msbm7
\font\tengoth=eufm10
\font\sevengoth=eufm7
\def\setA{{\bboard A}}
\def\setAA{\hbox{\sevenbboard A}}
\def\setZ{{\bboard Z}}

\def\calB{{\cal B}}
\def\calA{{\cal A}}
\def\calR{{\cal R}}
\def\calI{{\cal I}}
\def\calM{{\cal M}}
\def\calF{{\cal F}}
\def\calQ{{\cal Q}}
\def\calH{{\cal H}}
\def\calU{{\cal U}}
\def\calE{{\cal E}}
\def\Sym{\hbox{\tengoth S}}
\def\SSym{\hbox{\sevengoth S}}
\def\modI#1{\!\!\!\pmod{{{\cal I}#1}}}
\def\Airr{{\cal A}_{\rm irr}}
\def\Birr{{\cal B}_{\rm irr}}

\def\len{{\ell}}
\def\inv{{\rm inv}} 
\def\den{{\rm den}}
\def\exc{{\rm exc}}
\def\Ferm{\mathop{\rm Ferm}\nolimits}
\def\Bos{\mathop{\rm Bos}\nolimits}

\def\bw#1#2{{#1\choose #2}}

\def\smallmatrix#1{\vcenter{\offinterlineskip
     \halign{\vrule height 4pt depth 2pt width 0pt
      \hfil$\scriptstyle##$\hfil&&\kern
       3pt\hfil$\scriptstyle##$\hfil \crcr#1\crcr}}}

\catcode`\@=11
\def\matrice#1{\null \,\vcenter {\normalbaselines \m@th
\ialign {\hfil $##$\hfil &&\  \hfil $##$\hfil\crcr
\mathstrut \crcr \noalign {\kern -\baselineskip } #1\crcr
\mathstrut \crcr \noalign {\kern -\baselineskip }}}\,}

\def\pmatrice#1{\left(\null\vcenter {\normalbaselines \m@th
\ialign {\hfil $##$\hfil &&\thinspace  \hfil $##$\hfil\crcr
\mathstrut \crcr \noalign {\kern -\baselineskip } #1\crcr
\mathstrut \crcr \noalign {\kern -\baselineskip }}}\right)}

\def\bmatrice#1{\left[\null\vcenter {\normalbaselines \m@th
\ialign {\hfil $##$\hfil &&\thinspace  \hfil $##$\hfil\crcr
\mathstrut \crcr \noalign {\kern -\baselineskip } #1\crcr
\mathstrut \crcr \noalign {\kern -\baselineskip }}}\right]}

\catcode`\@=12
\def\modegale#1{\smash{\;\displaystyle\mathop{\equiv}^#1\;}}



%

\def\refBe{Be78}
\def\refCa{Ca72}
\def\refCl{Cl97}
\def\refCF{CF69}
\def\refFo{Fo65}
\def\refFHa{FH05a}
\def\refFHb{FH05b}
\def\refGLZ{GLZ05}
\def\refHa{Ha92}
\def\refMa{Ma15}
\def\refZ{Z80}
\def\refPW{PW91}

\rightline{2005/12/07 15:50}

\bigskip
\bigskip
\bigskip

\centerline{\bf Specializations and Extensions}
\centerline{\bf of the quantum MacMahon Master Theorem}
\medskip
\bigskip
\centerline{Dominique Foata and Guo-Niu Han}
\bigskip
\abstract{
We study some specializations and extensions of the quantum
version of the MacMahon Master Theorem derived by Garoufalidis,
L\^e and Zeilberger.  In particular, we obtain a $(t,q)$-analogue
for the Cartier-Foata noncommutative version and a  semi-strong
$(t,q)$-analogue for the contextual algebra.}
\medskip
\resume{
Nous \'etudions certaines sp\'ecialisations et extensions de la
version quantique du ``Master Theorem" de MacMahon,
\'etabli par Garoufalidis, L\^e et Zeilberger.
En particulier, nous obtenons un $(t,q)$-analogue
pour la version non-commutative de Cartier-Foata et un
$(t,q)$-analogue semi-fort pour l'lag\`ebre
contextuelle.
}

\bigskip
\section{1. Introduction} 
The Master Theorem derived by MacMahon [\refMa, vol.~1, p.~97] is a
fundamental result in Combinatorial Analysis and has had many
applications. To state it in a context relevant to its further
noncommutative extensions we use the following notations. Let~$r$
be a positive integer and~$\setA$  the {\it alphabet}
$\{1,2,\ldots,r\}$. A {\it biword} on~$\setA$  is a $2\times n$
matrix~$\alpha={x_1\cdots x_n\choose a_1\cdots a_n}$
($n\geq 0$), whose entries are in $\setA$, the first (resp. second)
row being called the {\it top word} (resp. {\it bottom word\/})
of~$\alpha$. The number $n$ is the {\it length} of
$\alpha$; we write $\len(\alpha)=n$. The biword $\alpha$
can also be viewed as a {\it word} of {\it biletters}
${x_1\choose a_1}\cdots{x_n\choose a_n}$, those biletters
$x_i\choose a_i$ being pairs of integers written vertically
with $x_i,a_i\in \setA$ for all $i=1,\ldots,n$. Let $\setB$ denote
the set of biletters and~$\calB$ the set of all biwords.

For each word $w=x_1x_2\ldots x_m\in \setA^m$
let $\inv\, w$ designate the number of inversions in~$w$, that is,
the number of pairs $(i,j)$ such that $1\le i<j\le m$ and
$x_i>x_j$. Also let $\overline w=x_{j_1}x_{j_2}\ldots x_{j_m}$ be
the {\it non-decreasing} rearrangement of~$w$. The $q$-{\it
Boson} is defined to be the infinite sum
$$
{\rm Bos}(q):=\sum_w  q^{\inv\,w}\pmatrice { \overline w\cr w\cr}
$$
over all words~$w$ from the free monoid ${\bboard A}^*$
generated by~$\bboard A$.
The $q$-{\it Fermion} is defined by
$$\Ferm(q)
:=\sum_{J\subset \setAA}
(-1)^{|J|} \sum_{\sigma\in {\SSym}_J}
(-1)^{\inv\, \sigma}
\pmatrice  {\sigma (i_1)&\sigma (i_2)&\cdots&\sigma (i_l)\cr
i_1&i_2&\cdots&i_l\cr},
$$
where $J=\{i_1<i_2<\cdots <i_l\}$ is a subset of~$\setA$ of
cardinality $|J|=\ell$ and~${\Sym}_J$ is the permutation group
acting on~$J$.

\goodbreak
When $q=1$ and when the biletters $x\choose a$ are assumed to
commute, the original MacMahon Master Theorem ({\it op. cit.})
asserts that the identity
$$
\Ferm(q)\times\Bos(q)=1\leqno\hbox{(\rm qMM)}
$$ 
holds. In the further noncommutative versions of the MacMahon
Master Theorem (for an arbitrary~$q$ or for $q=1$) specific
commutation rules for the biletters $x\choose a$ must be set. For
instance, in the Cartier-Foata version [\refCF] the biletters
$x\choose a$ and $y\choose b$ commute only if $x\not =y$ and
with $q=1$ the above relation holds in the large algebra of the
so-called monoid generated by those commutation rules.
Recently, Garoufalidis, L\^e and Zeilberger have established another
noncommutative version [\refGLZ] (called ``quantum Master
Theorem") by using  difference-operator techniques developed by
Zeilberger [\refZ]. They have borrowed their commutation rules
(see the reduction system $(SR_q)$ below) from the classical theory
of Quantum Groups. The (qMM)-identity holds in a quotient
algebra derived by those commutation rules, called the
{\it $q$-right quantum algebra}. The identity also holds in another
quotient algebra called the 1-{\it right quantum algebra}. Those
two identities were reproved in our previous paper [\refFHb] using a
different approach.

As other algebras previously introduced appear to be subalgebras of
the $q$-right or 1-right quantum algebras, it is
natural to determine whether the (qMM)-identity still hold in those
algebras and in which form. We can then show that there is a true
$(t,q)$-analogue of the Cartier-Foata noncommutative version
(\S\kern2pt6) and also a (qMM)-identity  in the
{\it contextual} algebra introduced in [\refHa] (\S\kern2pt 7).

The two right quantum algebras are defined by means of a {\it
reduction system} $(S)$ (see \S\S\kern2pt 2,~3), which
make it possible to construct an explicit {\it leftmost reduction}
$\alpha\mapsto [\alpha]_{S}$
mapping each biword~$\alpha$ onto a linear combination of
so-called {\it irreducible} biwords. We say that there is a {\it
strong} quantum Master Theorem, if the further identity
$[\Ferm(q)\times \Bos(q)]_S=1$ holds. This means that if
$(F\!B)_n$ denotes the linear combination of all biwords of
length~$n$ in the expansion of the product
$\Ferm(q)\times \Bos(q)$, then $(F\!B)_0=1$ and
for every $n\ge 1$ the expression~$(F\!B)_n$ can be
reduced to zero by applying the leftmost reduction inductively.

Borrowing the theory from Bergman [\refBe] we give the
definition of a reduction system in section~2. In section~3 we
recall the results on the (qMM)-identity that holds in the right
quantum algebras in its weak and strong forms. The
subalgebras of the right quantum algebra are introduced in
section~4. The list of results is presented in a table in section~5.
The $(t,q)$-analogue of the classical Master Theorem is proved in
section~6. Finally, we discuss in section~7 what we mean by
semi-strong analogue of the Master Theorem for the contextual
algebra.


\section{2. Reduction system and quotient algebra} 

Let $\setZ$ be the ring of all integers and $\setZq$ the ring of
polynomials in the variables~$q$, $q^{-1}$ submitted to
the rule $qq^{-1}=1$ with integral coefficients. The set
$\calA=\setZ\langle\!\langle \calB\rangle\!\rangle$
of the formal sums $\sum_{\alpha} c(\alpha) \alpha$, where
$\alpha\in\calB$ and $c(\alpha)\in\setZ$, together with the above
biword multiplication, the free addition and the free scalar product,
forms an algebra over $\setZ$, called the {\it free biword
large $\setZ$-algebra}. Similarly,
let $\calA_q=\setZq\langle\!\langle\calB\rangle\!\rangle$
denote the large $\setZq$-algebra of the
formal sums $\sum_\alpha c(\alpha)\alpha$ with
$c(\alpha)\in \setZq$ for all $\alpha\in {\cal B}$.
The subset of the {\it finite} formal sums is a subalgebra denoted
by 
$\setZq\langle\calB\rangle$.
The formal sums
$\sum_\alpha c(\alpha)\alpha$ in all those algebras
are called {\it expressions}.

Following Bergman's method [Be78]
each finite set of pairs
$(\alpha,E_\alpha)\in {\cal B}\times \setZq\langle\calB\rangle$
such that the {\it first components} $\alpha$ are distinct
is called a {\it reduction system}.
Let $F(S)=\{\alpha \mid (\alpha, E_\alpha)\in S\}$.
With $S$ we can associate an {\it oriented graph} $G=(\calB, \calE)$
defined as follows. The set of {\it vertices} is the set $\calB$
of all biwords. There is an oriented edge from $\beta$ to $\beta'$
if the following two conditions hold:

(1) there is a factorization
$\beta=\beta_1\alpha\beta_2$ such that $\beta_1,\alpha,\beta_2\in\calB$
with $\alpha\in F(S)$;

(2) $\beta'$ occurs in the expansion of the expression
$\beta_1 E_\alpha \beta_2$ with a non zero coefficient.

We assume that all positively oriented paths in the graph $G$ terminate.
This assumption, usually called the
{\it descending chain condition} (see [\refBe]), is assumed to hold
for all reduction systems  defined in the sequel.

A biword $\beta$ is said to be {\it irreducible}
if each factor $\alpha$ of $\beta$ is not a first component
of $S$. The set of irreducible biwords is denoted by $\Birr(S)$.
An expression $\sum_\alpha c(\alpha) \alpha$ is said to be
{\it irreducible} if $c(\alpha)=0$ for all $\alpha\not\in\calB_{\rm
irr}$. The set of irreducible expressions is denoted by
${\cal A}_{\rm irr}(S)$. Associated with the reduction system $S$ a
linear mapping $[\ ]_S: E\mapsto [E]_S$
of~$\calA_q$ onto the $\Airr(S)$, called the
{\it leftmost-reduction}, is defined by the following axioms:

\smallskip
(C1) The leftmost-reduction is $\setZq$-linear, i.e. for
$E_1, E_2\in \calA_q$ and $c_1, c_2\in\setZq$ we have
$[c_1 E_1+c _2 E_2]_S=c_1 [E_1]_S+c_ 2[E_2]_S$;
\smallskip

(C2) For every $\beta\in\Birr(S)$ we have $[\beta]_S=\beta$;
\smallskip

(C3) For every $(\alpha, E_\alpha)\in S$ we have $[\alpha]_S=E_\alpha$;
\smallskip

(C4) Let $\beta\not\in \Birr$ be a reducible biword,
so that $\beta$ can be factorized in a unique manner as
$\beta=\beta_1\beta_2{x\choose a}\beta_3$,
where $\beta_1, \beta_2, \beta_3\in\calB$, ${x\choose a}\in \setB$,
$\beta_1\beta_2\in\Birr$ and $\beta_2{x\choose a}\in F(S)$.
Then $[\beta]_S=[\ \beta_1 \ [\beta_2{x\choose a}]_S \ \beta_3]_S$.

\smallskip
Although the leftmost-reduction is recursively defined, it is
well defined. Starting with a biword $\beta$ and applying (C3)
finitely many times an irreducible expression is derived
(this is the {\it descending chain condition}, true by assumption).
When the biword~$\beta$ has more than one factor $\alpha\in F(S)$,
condition (C4) says that condition (C3) must be applied at the occurrence of
the {\it leftmost} factor $\alpha$. Consequently, the final
irreducible expression is unique.

\smallskip
Some reduction systems in the sequel will satisfy the further
condition:

\smallskip
(C5) Let $\beta\in \calB$ be a biword and
$\beta=\beta_1\beta_2\beta_3$ be any factorization
where $\beta_1, \beta_2, \beta_3\in\calB$.
Then $[\beta]_S=[\ \beta_1 \ [\beta_2]_S \ \beta_3]_S$.

\smallskip
A reduction system $S$ is said to be {\it reduction-unique}
if it satisfies condition (C5).
\smallskip
In practice, a reduction system is written as the set of equations
(or {\it commutation relations}) $\{\alpha=E_\alpha\}$.
Let $\calI(S)$  be the two-sided ideal of~$\calA_q$ generated by the
elements $(\alpha-E_\alpha)$  such that $(\alpha,E_\alpha)\in S$. Our main
algebraic structure is the {\it quotient algebra}
$\calU(S)={\calA_q}/\calI(S)$,
that will be studied for several reduction systems.


\section{3. Master Theorems on the right quantum algebra} 
\medskip
The reduction system $SR$ is defined by
$$ \cases{
{xy\choose aa} = {yx\choose aa}, & ($x>y$, all $a$);\cr
\noalign{\smallskip}
{xy\choose ab} = {yx\choose ba} + {yx\choose ab}-{xy\choose ba},
&$(x>y, a>b)$. \cr
}\leqno{(SR)}
$$
The quotient algebra $\calR:=\calU(SR)=\calA_q/\calI(SR)$ is called the
{\it $1$-right quantum algebra}.

\proclaim Theorem 1 {\rm ($1$-quantum Master Theorem)}.
The following identity holds:
$$
\Ferm(1) \times \Bos(1)\equiv 1 \modI{(SR)}.
$$

The above Theorem was first proved in [\refGLZ]. We can find another proof
in [\refFHb].  
\medskip
\proclaim Theorem 2 {\rm (Strong $1$-quantum Master Theorem)}.
The following identity holds:
$$
[\Ferm(1) \times \Bos(1)]_{SR} =1,
$$
where $[\ ]_{SR}$ is the leftmost-reduction associated with the reduction
system $(SR)$.

This Theorem was proved in [\refFHa] and is based on the following
property.

\proclaim Theorem 3 {\rm (Reduction-uniqueness)}.
The reduction system $SR$ is re\-duction-unique.

The above definitions and results have the following $q$-analogues.
Consider the reduction system
$$ \cases{
{xy\choose aa} = q{yx\choose aa}, & ($x>y$, all $a$);\cr
\noalign{\smallskip}
{xy\choose ab} = {yx\choose ba} + q{yx\choose ab}-q^{-1}{xy\choose ba},
&$(x>y, a>b)$. \cr
}\leqno{(SR_q)}
$$
The quotient algebra $\calR_q:=\calU(SR_q)=\calA_q/\calI(SR_q)$ is called
the
{\it $q$-right quantum algebra}.

\proclaim Theorem $1q$ {\rm ($q$-quantum Master
Theorem)}. We have the following identity:
$$
\Ferm(q) \times \Bos(q) \equiv 1 \modI{(SR_q)}.
$$

This Theorem was first proved in [\refGLZ]. We can find another proof
in [\refFHa] by using the ``$1=q$" principle.
For each biword $\bw uv$ define the statistic ``$\inv^-$" by
$$
\inv^- \bw uv =\inv (v)-\inv (u).
$$
Notice that ``$\inv^-$" may be negative.
The ``$1=q$" principle is based on the
weight function~$\phi$ defined for each biword $\alpha$
by $\phi(\alpha)=q^{\inv^-(\alpha)}\alpha$ and extended
to all of~${\cal A}_q$ by linearity.
A {\it circuit} is defined to be a biword whose top word
is a rearrangement of the letters of its bottom word.
An expression
$E=\sum_\alpha c(\alpha)\alpha$ is said to be {\it circular}
if  $c(\alpha) =0$ except when $\alpha$ is circular.
\proclaim Theorem 4 {\rm (``$1=q$" principle)}.
We have
$$
\leqalignno{
\phi([E]_{SR}) &= [\phi(E)]_{SR_q} \hbox{\quad for $E\in\calA_q$}; &\cr
\phi(E F) &= \phi(E) \phi(F) \hbox{\quad if $E$ and $F$ are two circular
expressions}. &\cr
}
$$

As proved in [\refFHa],
Theorems 3 and $1q$  imply the following theorem.

\proclaim Theorem $2q$ {\rm (Strong $q$-quantum Master
Theorem)}. We have
$$
[\Ferm(q) \times \Bos(q)]_{SR_q} =1,
$$
where $[\ ]_{SR_q}$ is the leftmost-reduction associated with the reduction
system $(SR_q)$.


\section{4. Subalgebras of the right quantum algebra}
Let $S$ be a reduction system
and let $\calU(S)=\calA_q/\calI(S)$ be the quotient algebra as
defined in section~2. If the commutation relations
$\alpha=E_\alpha$ derived from the reduction system $SR$ (resp.
$SR_q$) hold whenever the commutation relations
$\alpha'=E'_{\alpha'}$ derived from~$S$ hold, then $\cal U(S)$ is
a subalgebra of~$\calR$ (resp. of~$\calR_q$).
We examine several subalgebras $\calU(S)$ of~$\calR$
and~$\calR_q$ derived by specific reduction systems~$S$. In most
cases the implication
$\{\alpha'=E'_{\alpha '}\}\Rightarrow \{\alpha=E_{\alpha }\}$ is
immediate.

\smallskip
(a) {\it The commutative algebra} $\calM=\calU(SM)$ is defined by
means of the reduction system
$$
\leqalignno{
 {xy\choose ab} &= {yx\choose ba} \quad \hbox{($x>y$, all $a,b$).}
&{(SM)}\cr
 }
$$
This means that all biletters commute.
The algebra~$\calM$ is a subalgebra of~$\calR$. The original
Master Theorem identity due to MacMahon [\refMa] holds in $\calM$
(see also [\refCa]).

\smallskip
(b) {\it The Cartier-Foata algebra} $\calF=\calU(SF)$ is defined
by means of the reduction system
$$
\leqalignno{
 {xy\choose ab} &= {yx\choose ba}\quad \hbox{($x>y$).}
&(SF)\cr
 }
$$
The algebra $\calF$ was first considered by Cartier
and Foata [\refCF]. It is a subalgebra of~$\cal R$. A version
of~$\calF$ was also used in [\refFo] to derive a noncommutative
version of the  MacMahon Master Theorem.

\smallskip
(c) {\it The Cartier-Foata $q$-algebra} $\calF_q=\calU(SF_q)$ is
a $q$-version of the previous algebra. It is defined by means of
the reduction system
$$
\cases{
 {xy\choose ab} = {yx\choose ba}, &{if $x>y$ and $a>b$};\cr
 \noalign{\smallskip}
 {xy\choose ab} = q{yx\choose ba}, &{if $x>y$ and $a=b$};\cr
 \noalign{\smallskip}
 {xy\choose ab} = q^2{yx\choose ba}, &{if $x>y$ and $a<b$}.\cr
 }
\leqno{(SF_q)}
$$
The algebra $\calF_q$ is a subalgebra of $\calR_q$.

\smallskip
(d) {\it The quantum algebra} $\calQ_q=\calU(SQ_q)$ can only be
defined in its
$q$-version (see, for example, [\refPW]).
The underlying reduction system is given by
$$ \cases{
 {xy\choose aa} = q{yx\choose aa},&if $x>y$; \cr
 \noalign{\smallskip}
 {xx\choose ab} = q{xx\choose ba},&if $a>b$;\cr
 \noalign{\smallskip}
 {xy\choose ab} = {yx\choose ba},&if $x>y$ and $a<b$;\cr
 \noalign{\smallskip}
{xy\choose ab} =
  {yx\choose ba} + (q-q^{-1}){yx\choose ab}, &if $x>y$ and
$a>b$. \cr }\leqno{(SQ_q)}
$$
Notice that $\calQ_q$ is a subalgebra of $\calR_q$ and
$\calQ_q|_{q=1} = \calM$.

\smallskip
(e) {\it The contextual algebra} $\calH=\calU(SH)$ is defined by the
reduction system [\refHa]
$$
\cases{
 {xy\choose ab} = {yx\choose ab},&if $x>y$ and $V(x,y,a,b)>0$;\cr
 \noalign{\smallskip}
 {xy\choose ab} = {yx\choose ba},&if $x>y$ and $V(x,y,a,b)<0$;\cr
 }\leqno{(SH)}
$$
where $V(x,y,a,b)=(a-x-1/2)(a-y-1/2)(b-x-1/2)(b-y-1/2)$.
In Proposition~6 we show that the relations displayed in $(SH)$
imply the relations~displayed in~$(SR)$. The
algebra~$\calH$ is then a subalgebra of~$\calR$. By convention, we
define the $q$-version of~$\calH$ as being $\calH_q=\calH$
itself. Therefore,~$\calH_q$ is not a subalgebra of $\calQ_q$.

\medskip

\section{5. List of results} 
In sections 2 and 3 we have stated six theorems:
Theorem~1---Theorem~4 and Theorem~$1q$, Theorem~$2q$. They all
relate to the algebras~$\calR$ or~$\calR_q$ (see the last row in
the following table). On the other hand,  three subalgebras~$\calM$,
$\calF$, $\calH$ of~$\calR$ and two subalgebras~$\calQ_q$,
$\calF_q$  of~$\calR_q$ have been introduced. Our purpose is to
state and prove, whenever possible, the analogues of those theorems
for those subalgebras. As can be seen in the table
twenty-six statements are to be made. Each of them refers to an
algebra (the row index) and a Theorem (the column index).

$${
\def\tvi{\vrule height 12pt depth 5pt width 0pt}
\def\tv{\tvi\vrule}

\def\oui{$\otimes$}
\def\demi{$\oslash$}
\def\non{$-$}
\vbox{\offinterlineskip\halign{
\ #\hfil\ \tv&&\strut\ #\hfil\ \cr
             &Th.1          &Th.2            &Th.3          &Th.4
&Th.$1q$              &Th.$2q$            &\tv&         \cr
\noalign{\hrule}
$\calM$       &\oui$^{(7)}$  &\oui$^{(15)}$   &\oui$^{(11)}$ &\non
&\non               &\non            &\tv&\non  \cr
$\calQ=\calM$ &\oui$^{(8)}$  &\oui$^{(16)}$   &\oui$^{(13)}$ &\non
&\oui$^{(17)}$      &\oui$^{(21)}$   &\tv&$\calQ_q$  \cr
$\calF$       &\oui$^{(9)}$  &\oui$^{(14)}$   &\oui$^{(12)}$ &\oui$^{(18)}$
&\oui$^{(19)}$      &\oui$^{(20)}$  &\tv&$\calF_q$  \cr
$\calH$       &\oui$^{(10)}$ &\non$^{(23)}$   &\non$^{(22)}$ &\oui$^{(24)}$
&\oui$^{(25)}$      &\demi$^{(26)}$ &\tv&$\calH_q=\calH$  \cr
$\calR$       &\oui$^{(1)}$  &\oui$^{(5)}$    &\oui$^{(3)}$  &\oui$^{(4)}$
&\oui$^{(2)}$       &\oui$^{(6)}$    &\tv&$\calR_q$  \cr
}}
}$$

Comments on the content of the table:

(1) This is Theorem 1 (quantum Master Theorem), first proved in
[\refGLZ] and reproved in [\refFHb].

(2) This $q$-version of the quantum Master Theorem for the right
quantum algebra can be deduced from (1) and the ``$1=q$"
principle (4) [\refFHa].

(3) See [\refFHa]. 

(4) See [\refFHa]. The weight function for the ``$1=q$" principle is
``$\inv^-$".

(5) This is Theorem 2 (strong 1-quantum Master Theorem). It can
be deduced from (1) and~(3) (see [\refFHa]).

(6) This $q$-version of the strong quantum Master Theorem for
the right quantum  algebra can be deduced from (2) and (3).

(7-10) All those theorems are corollaries of (1) because $\calM$,
$\calQ$, $\calF$, $\calH$ are all subalgebras of $\calR$.
Also, notice that $(7)=(8)$.

(11-12) The bases of those two algebras are easy to construct

(13) See Parshall and Wang [\refPW] in its $q$-version.
If $q=1$, then $(13)=(11)$.

(14) See Foata [\refFo] and Cartier-Foata [\refCF]. It can be deduced
from (9) and (12).

(15-16) As $q=1$ for those two cases, we have $(15)=(16)$. The
result can be deduced from (7) and (11).

(17) This is a corollary of (2). Notice that it {\it cannot} be
deduced from (8) because of the lack of any ``$1=q$" principle for
$\calQ$. See [\refGLZ].

(18) See section 6 in this paper. The weight function is ``$\inv^-$".

(19) This is a corollary of (2). It can also be deduced from (9)
and the ``$1=q$" principle for $\calQ$ (18).

(20) This $q$-version of the strong quantum Master Theorem for
the right quantum  algebra can be deduced from (19) and (12). See
section 6.

(21) This $q$-version of the strong quantum Master Theorem for
the right quantum  algebra can be deduced from (17) and (13).
See [\refGLZ].

(22-23) The {\it leftmost-reduction} process is not unique.
See section 7.
There is no {\it strong} Master Theorem for $\calH$.

(24)  The ``$1=q$" principle uses the weight function defined by the Denert
statistics [\refHa].

(25) It can be deduced from (10) and (24). But it is {\it not} a
corollary of (2), because $\calH_q$ is not a subalgebra of
$\calR_q$.

(26) Even though there is no strong Master Theorem for $q=1$ (23),
there is a {\it semi-strong} Master Theorem in its $q$-version.
See section 7.


\section{6. A $(t,q)$-version of the classical Master Theorem} 
In this section we derive a $(t,q)$-analogue of the MacMahon
Master Theorem for the classical noncommutative case [\refCF].
We could give a short proof based on the method used in [\refFHb],
because the further properties needed  are straightforward.
However we prefer to deduce this $(t,q)$-analogue identity from
the classical noncommutative version [\refCF]. To that end we
consider the reduction system $(SF_q)$ and the quotient algebra
$\calF_q$ defined in section 4. It is easy to see that $(SF_q)$ is
a {\it reduction-unique} system and that $\calF_q$ is a
subalgebra of~$\calR_q$.

The {\it $(t,q)$-Boson} is defined to be the infinite sum
$$
{\rm Bos}(t,q):=\sum_w t^{\exc {\overline w\choose w}}
q^{\inv\, w} \pmatrice { \overline w\cr w\cr}
$$
over all words~$w$ from the free monoid ${\bboard A}^*$
generated by~$\bboard A$.
The definition of the statistic ``$\exc$" (``number of
exceedances") for biwords is classical and can be found in [\refHa].
The {\it
$(t,q)$-Fermion} is an extension of the $q$-Fermion and
reads
$$\Ferm(t,q)
:=\sum_{J\subset \setAA}
(-1)^{|J|} \sum_{\sigma\in {\SSym}_J}
t^{\exc\,\alpha}(-q)^{\inv\, \sigma}
\pmatrice  {\sigma (i_1)&\sigma (i_2)&\cdots&\sigma (i_l)\cr
i_1&i_2&\cdots&i_l\cr}.
$$ 
In the above definition $\alpha$ denotes the biword
\smash{$\left(\smallmatrix  {\sigma (i_1)&\sigma
(i_2)&\cdots&\sigma (i_l)\cr i_1&i_2&\cdots&i_l\cr}\right)$}.

\smallskip
\proclaim Theorem 1F$q$ {\rm ($q$-quantum Master Theorem)}.
The following identity\hfil\break holds:
$$
\Ferm(t,q) \times \Bos(t,q) \equiv 1 \modI{(SF_q)}.
$$

\goodbreak

\proclaim Theorem 2F$q$ {\rm (Strong $q$-quantum Master
Theorem for $SF_q$)}. The following identity holds:
$$
[\Ferm(t,q) \times \Bos(t,q)]_{SF_q} =1,
$$
where $[\ ]_{SF_q}$ is the leftmost-reduction associated with the
reduction system $(SF_q)$.

Theorem 1F$q$ is an immediate consequence of Theorem 2F$q$. To
prove Theorem~2F$q$ we first recall the classical noncommutative
version of the Master Theorem [\refCF] that we express as follows.

\proclaim Theorem 2F {\rm (Classical noncommutative version)}.
The following identity holds:
$$
[\Ferm(1) \times \Bos(1)]_{SF} =1,
$$
where $[\ ]_{SF}$ is the leftmost-reduction associated with the
reduction system $(SF)$.

Now, introduce the weight function
$\phi_{SF}: \calA_q
\rightarrow \setZtq \langle\!\langle\calB\rangle\!\rangle$ defined
for each biword $\alpha$
by $\phi_{SF}(\alpha) = t^{\exc(\alpha)}q^{\inv^-(\alpha)}\alpha$
and extended to all of~${\cal A}_q$ by linearity.
The ``$1=q$" principle for $\calF$ is stated next as Theorem~4F.
The proof is straightforward and omitted.

\proclaim Theorem 4F {\rm (``$1=q$" principle for $\cal F$)}.
We have
$$
\leqalignno{
\phi_{SF}([E]_{SF}) &= [\phi_{SF}(E)]_{SF_q}
\hbox{\ \qquad for $E\in\calA_q$}; &\cr
\phi_{SF}(E F) &= \phi_{SF}(E) \phi_{SF}(F)
\hbox{\quad if $E$ and $F$ are two circular expressions.} &\cr
}
$$

We complete the proof of Theorem 2Fq as follows.
We have
$$\eqalignno{
[\Ferm(1) \times \Bos(1)]_{SF} &=1, &\hbox{[by Theorem 2F]}\cr
\phi_{SF}([\Ferm(1) \times \Bos(1)]_{SF}) &=\phi_{SF}(1)=1, &\cr
[\phi_{SF}(\Ferm(1) \times \Bos(1))]_{SF_q} &=1, &\hbox{[by Theorem 4F
(i)]}\cr
[\phi_{SF}(\Ferm(1)) \times \phi_{SF}(\Bos(1))]_{SF_q} &=1.
  &\hbox{[by Theorem 4F (ii)]}\cr
}
$$
Finally, it is easy to verify that $\phi_{SF}(\Ferm(1))=\Ferm(t,q)$
and also
$\phi_{SF}(\Bos(1))=\Bos(t,q)$. \qed


\section{7. Semi-strong $(t,q)$-analogue of the Master Theorem for
$\calH$} 
The contextual algebra~$\cal H$ introduced in
\S\kern2pt 4(e) provides an example~of~an algebra in which the
``qMM" theorem holds, but not the ``strong qMM" one. Those two
assertions are proved in this section. We also~obtain a
$(t,q)$-version of the quantum Master Theorem, that can be
regarded as a {\it semi-strong} extension, because the product
$\Ferm(1)\times
\Bos(1)$ is reduced to~1 after {\it two} steps:
 first, by a leftmost reduction (this is the strong~part),
second, by taking a homomorphic image of the product
thereby  reduced.

Consider the reduction system
$(SH)$ defined in section~4, a
system that has been extensively studied in [\refHa] and [\refCl].

\goodbreak
\proclaim Proposition 5.
The algebra $\calH$ is a subalgebra of $\calR$.

{\it Proof}.
If $x>y$ and $V(x,y,a,b)>0$, then
 $\displaystyle{xy\choose ab} = {yx\choose ab}$ by definition.
As $V(x,y,b,a)=V(x,y,a,b)>0$, we also have
 $\displaystyle{xy\choose ba} = {yx\choose ba}$.
On the other hand, if $x>y$ and $V(x,y,a,b)<0$, then
$\displaystyle{xy\choose ab} = {yx\choose ba}$.
As $V(x,y,b,a)=V(x,y,a,b)<0$, we also have
$\displaystyle{xy\choose ba} = {yx\choose ab}$.
In both cases 
$$
{xy\choose ab} = {yx\choose ba} + {yx\choose ab}-{xy\choose ba}.\qed
$$

\proclaim Proposition 6.
The reduction system $(SH)$ is not reduction-unique.

{\it Proof}. It suffices to provide this following counter-example:
$$
\leqalignno{
{321\choose213}\modegale1
{231\choose213}\modegale2
{213\choose231}\modegale1
{123\choose321} \modI{(SH)};&\cr
{321\choose213}\modegale2
{312\choose213}\modegale1
{132\choose123}\modegale2
{123\choose132} \modI{(SH)};&\cr
}
$$
where $\modegale{i}$ means that the reduction relations (SH) are
applied at positions $(i,i+1)$.
Both ${123\choose321}$ and ${123\choose132}$ belong to
$\Birr(SH)$ and are derived by reduction from the same circuit.\qed

\medskip
The integral-valued statistic ``den" is the {\it Denert statistic} for
biwords. Its definition, as well as the definition of ``exc," can be
found  in [\refHa].

\proclaim Proposition 7.
Let $u$ and $v$ be two circuits. If
$\alpha\equiv\beta\modI{(SH)}$, then
$$(\exc, \den)\alpha =(\exc, \den)\beta.$$

See [\refHa] for the proof. Note that Clark [\refCl] has shown that
the converse is also true for each pair of bipermutations $\alpha$,
$\beta$. Using Proposition 5 and Theorem 1 we can deduce  the
following ``1-quantum Master Theorem" for $\calH$.

\proclaim Theorem 1H {\rm (1-quantum Master Theorem for
$\calH$)}. The following identity holds:
$$
\Ferm(1) \times \Bos(1)\equiv 1 \modI{(SH)}.
$$

However there is no ``strong $q$-quantum Master Theorem" for
$\calH$, as shown in the next proposition.

\proclaim Proposition 8.
We have
$[\Ferm(1) \times \Bos(1)]_{SH} \not=1$,
where $[\ ]_{SH}$ is the leftmost-reduction associated with the
reduction system $(SH)$.

{\it Proof}. For $r=3$ we calculate
$$
\leqalignno{
[\Ferm(1) \times \Bos(1)]_{SH}
&=1
+{123\choose213}
-{123\choose312}
-{123\choose321}
+{123\choose132}\cr
&\qquad
+{1123\choose2113}
-{1123\choose3121}
-{1123\choose3112}
+{1123\choose1312}\cr
&\qquad
+{1223\choose2123}
-{1223\choose3122}
-{1223\choose3212}
+{1223\choose1322}\cr
&\qquad
+{1233\choose2133}
-{1233\choose3123}
-{1233\choose3213}
+{1233\choose1323}\cr
&\qquad+\cdots\cr
&\not=1.\qed\cr
}
$$

Now, let the weight function
$\psi: \calA_q \rightarrow \setZtq
\langle\!\langle\setA^*\rangle\!\rangle$ be defined for each biword
$\alpha={u\choose v}$ by 
$$\psi{u\choose v}=t^{\exc (\alpha)}q^{\den(\alpha)}u$$
and be extended to all of~${\cal A}_q$ by linearity.
From Theorem~1H and Proposition~7 we deduce the following
result.

\proclaim Theorem 2H {\rm (Semi-strong quantum Master Theorem
for $\calH$)}. We have
$$
\psi([\Ferm(1) \times \Bos(1)]_{SH}) =1,
$$
where $[\ ]_{SH}$ is the leftmost-reduction associated with the
reduction system $(SH)$.

For example, take $r=3$, we have
$$
\leqalignno{
\psi([\Ferm(1) \times \Bos(1)]_{SH})
&=1
+tq(123)
-tq(123)
-tq^2(123)
+tq^2(123)\cr
&\quad
+tq(1123)
-tq^2(1123)
-tq(1123)
+tq^2(1123)\cr
&\quad
+tq(1223)
-tq(1223)
-tq^2(1223)
+tq^2(1223)\cr
&\quad
+tq(1233)
-tq(1233)
-tq^2(1233)
+tq^2(1233)\cr
&\quad+\cdots\cr
&=1.\qed\cr
}
$$

\hfill\eject
\vglue 2mm
\bigskip

{
\eightpoint

\centerline{\bf References} 
\bigskip

\article \refBe|George M. Bergman|The Diamond Lemma for Ring
Theory|Adv. Math.|29|1978|178--218|
\smallskip

\livre \refCa|P. Cartier|
La s\'erie g\'en\'eratrice exponentielle, Applications probabilistes
et alg\'ebriques|
Publ. I.R.M.A., Universit\'e Louis Pasteur, Strasbourg,
{\oldstyle 1972}, 68 pages|
\smallskip

\article \refCl|Robert. J. Clarke|Han's conjecture on permutations|European
J. 
Combin.|18|1997|511--524|
\smallskip

\livre \refCF|P. Cartier, D. Foata|Probl\`emes combinatoires de
permutations et r\'earran\-ge\-ments|Berlin, Springer-Verlag,
{\oldstyle 1969} ({\sl Lecture Notes in Math., {\bf 85}})|
\smallskip

\article \refFo|D. Foata|Etude alg\'ebrique de certains probl\`emes
d'Analyse Combinatoire et du Calcul des Probabilit\'es|Publ.
Inst. Statist. Univ. Paris|14|1965|81--241|
\smallskip

\divers \refFHa|D. Foata, G.-N. Han|{\sl A basis for the right quantum
algebra
and the ``$1=q$" principle}, preprint, {\oldstyle 2005}|
\smallskip

\divers \refFHb|D. Foata, G.-N. Han|{\sl A New Proof of the
Garoufalidis-L\^e-Zeilberger Quantum MacMahon Master
Theorem}, preprint, {\oldstyle 2005}|

\divers \refGLZ|S. Garoufalidis, T. Tq L\^e, D. Zeilberger|
{\sl The Quantum MacMahon Master Theorem}, to be appear
in {\sl Proc. Natl. Acad. Sci.},
{\oldstyle 2005}|
\smallskip

\article \refHa|Guo-Niu Han|Une transformation fondamentale
sur les r\'earrangements de mots|Adv.
Math.|105|1994|26--41|
\smallskip
 
\livre \refMa|P. A. MacMahon|Combinatory Analysis {\rm vol.
I, II}|Cambridge Univ. Press, {\oldstyle 1915}. Reprinted by
Chelsea, New York, {\oldstyle 1960}|
\smallskip

\livre \refPW|B. Parshall, J.-P. Wang|Quantum linear groups|
Memoirs Amer. Math. Soc, {\bf 89}, {\oldstyle 1991}|
\smallskip

\article \refZ|D. Zeilberger|The algebra of linear partial difference
operators and its applications|SIAM J. Math. Anal.|11, no. 6|1980|919--932|
\smallskip


}

\bigskip\bigskip
\hbox{\vtop{\halign{#\hfil\cr
Dominique Foata \cr
Institut Lothaire\cr
1, rue Murner\cr
F-67000 Strasbourg, France\cr
\noalign{\smallskip}
{\tt foata@math.u-strasbg.fr}\cr}}
\qquad
\vtop{\halign{#\hfil\cr
Guo-Niu Han\cr
I.R.M.A. UMR 7501\cr
Universit\'e Louis Pasteur et CNRS\cr
7, rue Ren\'e-Descartes\cr
F-67084 Strasbourg, France\cr
\noalign{\smallskip}
{\tt guoniu@math.u-strasbg.fr}\cr}}
}

\vfill\eject

\bye